\documentclass[english]{smfart}
\usepackage{sabbah_rigid}
\begin{document}
\frontmatter
\title[Rigid irreducible meromorphic connections in dimension one]{Remarks on rigid irreducible meromorphic connections on the projective line}

\author[C.~Sabbah]{Claude Sabbah}
\address{CMLS, CNRS, École polytechnique, Institut Polytechnique de Paris, 91128 Palaiseau cedex, France}
\email{Claude.Sabbah@polytechnique.edu}
\urladdr{https://perso.pages.math.cnrs.fr/users/claude.sabbah}


\begin{abstract}
We illustrate the Arinkin-Deligne-Katz algorithm for rigid irreducible meromorphic bundles with connection on the projective line by giving motivicity consequences similar to those given by Katz for rigid local systems \cite{Katz96}.
\end{abstract}

\maketitle
\vspace*{-\baselineskip}
\tableofcontents
\mainmatter

\section{Introduction}
Let $k:U\hto\PP^1$ be the inclusion of a proper Zariski open subset of the complex projective line $\PP^1$ and let $(\cV,\nabla)$ be an irreducible algebraic bundle\footnote{In the following, bundle means algebraic vector bundle.} of rank $r$ with connection on $U$. We say that $(\cV,\nabla)$ is \emph{rigid} if any other $(\cV',\nabla)$ on $U$ having at each puncture $x\in\PP^1\moins U$ a~formal structure isomorphic to that of $(\cV,\nabla)$ satisfies $(\cV',\nabla)\simeq(\cV,\nabla)$. It is proved (\cite{Katz96,B-E04}) that, on $\PP^1$, this property is equivalent to \emph{cohomological rigidity}, \ie the \emph{rigidity index}
\[
\rig(V,\nabla):=2-h^1_{\dR}(\PP^1,k_{\dag+}\cEnd(V,\nabla))
\]
is equal to $2$, where $k_{\dag+}$ denotes the minimal extension in the sense of $\cD$-modules.

The Arinkin-Deligne-Katz algorithm \cite{Arinkin08}, which relies on the property that the rigidity index is preserved by Fourier transformation (\cite{B-E04}) provides an inductive way of checking whether a given irreducible $(\cV,\nabla)$ is rigid by means of successive specific transformations: $(\cV,\nabla)$ is rigid if and only if the sequence of transformations in the algorithm reaches the trivial rank-one bundle with connection ($\cO_{U'},\rd)$ on some open subset $U'\subset\PP^1$. On the other hand, there is a one-to-one correspondence \hbox{between} irreducible bundles with connection $(\cV,\nabla)$ on some Zariski open subset of~$\PP^1$ and irreducible holonomic $\cD$\nobreakdash-modules $\cM$ on $\PP^1$ by the inverse functors ``middle extension'' and ``restriction to a suitable Zariski open set'', and the algorithm works with the latter objects.

For $N\in\NN^*$, we say that a bundle with connection $(\cV,\nabla)$ (or its middle extension~$\cM$ to~$\PP^1$) is \emph{$N$-quasi-unipotent} if the eigenvalues of the formal monodromy at each \hbox{$x\in\PP^1\moins U$} belong to $\mu_N$ ($N$-th roots of the unity).

The results of this note concern quasi-unipotent rigid irreducible bundles with connection on some proper Zariski open subset $U\subsetneq\PP^1$. They consist of applications of the Arinkin-Deligne-Katz algorithm. A first application has already been given in \cite{Bibi15}, where it is shown that any rigid irreducible $(\cV,\nabla)$ on $U$ (without the condition of quasi-unipotency) underlies an irregular mixed Hodge module structure which is pure (of~weight zero, say).

The motivation for this question came from various recent talks by Michael Groechenig, Aaron Landesman and Daniel Litt on their respective works \cite{E-G17} and \cite{L-L22}. Of course, the technique used here on $\PP^1$ does not extend to higher dimensions, but it opens the way to questions in higher dimensions in the setting of irregular singularities. In another direction, the finiteness result of Haoyu Hu and Jean-Baptiste Teyssier \cite{H-T22} looks promising.

\subsubsection*{Rank one}
Any rank-one bundle with connection $(\cV,\nabla)$ on $U$ is isomorphic to one of the form \hbox{$(\cO_U,\nabla^\reg+\rd\varphi)$}, where $\nabla^\reg$ is a connection having regular singularities on $\PP^1\moins U$ and $\varphi$ is a regular function on $U$. It is clearly irreducible and is cohomologically rigid (because $\cEnd(\cV,\nabla)=(\cO_U,\rd)$). That it is rigid is seen as follows (a special case of the criterion mentioned above): if $(\cV'',\nabla):=(\cV,\nabla)^\vee\otimes(\cV',\nabla)$ is a rank-one local system which has regular singularity at each $x\in\PP^1\moins U$ and has trivial local monodromy there, it extends to a trivial bundle with connection on $\PP^1$.

The bundle with connection $(\cV,\nabla)$ is quasi-unipotent if and only if $(\cO_U,\nabla^\reg)$ is~so, and this amounts to the property that a suitable tensor power $(\cO_U,\nabla^\reg)^{\otimes N}$ is isomorphic to $(\cO_U,\rd)$.

\Subsubsection*{Finiteness}
\begin{definition}[Minimally ramified polar part]
Let $p$ be an integer $\geq1$ and let $\varphi\in\CC\lpr t^{1/p}\rpr/\CC\lcr t^{1/p}\rcr$ be a nonzero ramified polar part of ramification order $p$. We~say that~$\varphi$ is \emph{minimally ramified} if it is not the pullback of a ramified polar part in $\CC\lpr t^{1/p'}\rpr/\CC\lcr t^{1/p'}\rcr$ with $p'< p$.
\end{definition}

If $\varphi$ is a minimally ramified polar part, it~yields a $p$-dimensional $\CC\lpr t\rpr$-vector space with connection, that we denote by $\El(\varphi)$, obtained as the pushforward by $t^{1/p}\mto t$ of \hbox{$(\CC\lpr t^{1/p}\rpr,\rd+\rd\varphi)$}. By~the Levelt-Turrittin theorem, any finite-dimensional $\CC\lpr t\rpr$\nobreakdash-vector space with \hbox{connection} can be written in a unique way as the direct sum of terms $\El(\varphi)\otimes R_\varphi$, where~$\varphi$ runs in a finite set of minimally ramified polar parts and $R_\varphi$ is a finite-dimensional $\CC\lpr t\rpr$-vector space with regular singular connection (\cf\cite{B-E04}, \cf also \cite{Bibi07a}). The minimally ramified polar parts entering the Levelt-Turrittin decomposition are called the \emph{exponential factors} (also called its \emph{irregular values}) of the $\CC\lpr t\rpr$-vector space with connection. Given a finite set $\Phi$ of (possibly non minimally) ramified polar parts, we say with some abuse that the exponential factors of $(\cV,\nabla)$ belong to $\Phi$ if any suitably further ramified exponential factor of $(\cV,\nabla)$ belongs to $\Phi$.

\begin{theoremealph}[Finiteness]\label{th:finiteness}
Given integers $r,N\geq1$ and a finite set $\Phi$ of ramified polar parts, there exists only a finite number of quasi-unipotent rigid irreducible bundles with connection $(\cV,\nabla)$ on $U$ of rank $r$ such that
\begin{itemize}
\item
the order of quasi-unipotency is at most $N$,
\item
At each $x\in\PP^1\moins U$, the exponential factors of (the formalization of) $(\cV,\nabla)$ at $x$ belong to $\Phi$.
\end{itemize}
\end{theoremealph}

\subsubsection*{Exponential-geometric origin}
Let $Z$ be a smooth complex quasi-projective variety. We say that an algebraic vector bundle with an integrable connection $(\cV,\nabla)$ on $Z$ is \emph{of exponential-geometric origin} if there exist a Zariski dense open subset \hbox{$j:U\hto Z$}, a~morphism $f:Y\to U$ from a smooth quasi-projective variety and a regular function~$\varphi$ on~$Y$ such that $j^*(\cV,\nabla)$, regarded as an holonomic $\cD_U$-module, is a submodule of $\cH^kf_+(\cO_Y^r,\rd+\rd\varphi)$ for some $r\geq1$ and some $k\in\ZZ$. (This is an adaptation of the definition of ``geometric origin'' in \cite{E-K21}, \cf also \cite{L-L22}).

Roughly speaking, horizontal sections (or solutions) of such a $(\cV,\nabla)$ on $U^\an$ can be given an integral expression, with the integrand being of the form $\rme^\varphi\cdot\omega$ for some algebraic differential form $\omega$ (\cf \cite{H-R08}).

\begin{theoremealph}\label{th:expmot}
Any quasi-unipotent rigid irreducible $(\cV,\nabla)$ on $U$ has expo\-nential-geometric origin.
\end{theoremealph}

\subsubsection*{Integral structure}
To any bundle with connection $(\cV,\nabla)$ on $U$ is associated a Stokes-filtered local system $(\cL,\cL_\bbullet)$ on the oriented real blow-up space $\wt\PP^1$ of $\PP^1$ at the punctures $\PP^1\moins U$ (\cf\cite{Malgrange91,Bibi10}). The local system $\cL$ and the terms $\cL_\bbullet$ of the Stokes filtration are sheaves of $\CC$\nobreakdash-vector spaces. We say that the $\CC$-Stokes-filtered local system has an integral structure if it comes by extension of scalars from $\QQ[\zeta]$ to $\CC$ (some $\zeta\in\mu_N$ and some $N\geq1$) from a $\QQ[\zeta]$-Stokes-filtered local system with a $\ZZ[\zeta]$-structure in the sense of Definition \ref{def:ostructure}.

\begin{theoremealph}[Integral structure]\label{th:integrality}
The Stokes-filtered local system associated to any quasi-unipotent rigid irreducible bundle with connection $(\cV,\nabla)$ on $U\subsetneq\PP^1$ has an integral structure.
\end{theoremealph}

\begin{remarque}
One can define the notions of irreducibility, rigidity and quasi-unipotency for a Stokes-filtered local system. Due to the Riemann-Hilbert correspondence of Deligne and Malgrange (\cf \cite{Malgrange91}), they~correspond to those of the associated bundle with connection. The previous proposition can be stated as the property that a quasi-unipotent rigid irreducible Stokes-filtered local system has an integral structure.
\end{remarque}

\begin{exemple}\label{ex:hypergeom}
A (possibly confluent) non resonant hypergeometric differential equation is irreducible and rigid, and it is quasi-unipotent if its local exponents belong to $\frac1N\ZZ$ for some $N\in\NN^*$. In \cite{D-M89} and \cite{Hien22}, the authors compute the Stokes matrices of confluent hypergeometric equations and the existence of an integral structure is then clear from their formulas. On the other hand, in \cite{B-H-H-S22}, as a particular case of their results, the authors make explicit its exponential-geometric origin and show that the associated enhanced ind-sheaf is defined on a cyclotomic extension of~$\QQ$.
\end{exemple}

\begin{exemple}\label{ex:G2}
In \cite{Jakob20}, the author classifies rigid irreducible bundles with connection $(\cV,\nabla)$ on $\Gm$ with an irregular singularity at infinity of slope one at infinity and differential Galois group $G_2$. The family he obtains depends on various parameters in $\CC^*$. Quasi-unipotency is equivalent to these parameters belonging to $\mu_N$ for some $N\geq1$, and Property \ref{th:integrality} claims that, in such a case, the corresponding Stokes-filtered local system has an integral structure. See also \cite[\S8.4]{Katz96} and \cite{D-R10} for the geometric origin in the tame case.
\end{exemple}

\section{Finiteness}
In this section we prove Property \ref{th:finiteness}. We consider the data $(U,N,\Phi)$ and quasi-unipotent rigid irreducible bundles with connection of rank $r$ with these data, \ie defined on~$U$, quasi-unipotent of order dividing $N$ and with exponential factors contained in~$\Phi$.

A first approach to Property \ref{th:finiteness} is by noticing that from the data $(U,N,\Phi)$, one can cook up only a finite number of possible formal structures $\bigoplus[\El(\varphi)\otimes R_\varphi]$ at each $x\in \PP^1\moins U$, with $\varphi\in\Phi$, $R_\varphi$ being regular and $N$-quasi-unipotent. For each such data of formal structures at every $x\in\PP^1\moins U$, there exists a smooth affine moduli space of finite type over $\CC$ such that the corresponding irreducible rigid $(V,\nabla)$ are isolated points of this space (\cf\cite[Proof of Th.\,4.10]{B-E04}). Finiteness follows.

We will now prove finiteness as a consequence of the Arinkin-Deligne-Katz algorithm. This methods, being more constructive, is more quantitative, although it uses the equivalence between rigidity and cohomological rigidity shown in \cite[Th.\,4.10]{B-E04}. The proof is by induction on the rank $r$ of $\cV$, and we denote by $\eqref{th:finiteness}_r$ the statement that Property \ref{th:finiteness} holds for bundles having any set of data $(U,N,\Phi)$ and of rank $r$.

\begin{proof}[Proof of $\eqref{th:finiteness}_1$]
We choose an affine coordinate $t$ on $\PP^1$ such that \hbox{$\infty\in U$}. For each $x\in \PP^1\moins U$, we regard $\varphi_x\in\Phi$ as a polynomial in $1/(t-x)$ with no constant term. Any choice of a family $(\varphi_x)_{x\in\PP^1\moins U}$ of elements of $\Phi$ (there are finitely many such choices) yields a unique regular function $\varphi$ on $U$ (namely, $\varphi(t)=\sum_{x\in\PP^1\moins U}\varphi_x(t-x)$). Given a bundle with connection $(\cV,\nabla)$ of rank one on~$U$ with data $N,\Phi$, there exists such a family $(\varphi_x)_{x\in\PP^1\moins U}$ such that $(\cV,\nabla-\rd\varphi)$ has regular singularities at each point $x\in\PP^1\moins U$. Since the local monodromies at each such~$x$ belong to $\mu_N$, there is only a finite number of such bundles with connection.
\end{proof}

We now assume $r\geq2$ and $\eqref{th:finiteness}_{<r}$, and we will prove $\eqref{th:finiteness}_r$. We are given $(U,\Phi,N)$ and we will prove $\eqref{th:finiteness}_r$ for these data, a statement that we denote by $\eqref{th:finiteness}_r(U,\Phi,N)$.

\begin{itemize}
\item
We can (and will) assume that $\#(\PP^1\moins U)\geq3$.
\end{itemize}
Indeed, if $\#(\PP^1\moins U)<3$, let $U_1\subset U$ with $\#(\PP^1\moins U_1)\geq3$. By considering $U\moins U_1$ as apparent singularities for $(\cV,\nabla)$, we have the implication
\[
\eqref{th:finiteness}_r(U_1,\Phi,N)\implies \eqref{th:finiteness}_r(U,\Phi,N).
\]

Given $(\cV,\nabla)$ rigid irreducible of rank $r$ on $U$, there exists at most one $x\in\PP^1\moins U$ where $(\cV,\nabla)_{\wh x}$ is special with respect to the Katz algorithm (Case II of \cite[\S4.1]{Arinkin08}). By possibly adding an apparent singularity, we can thus assume that
\begin{itemize}
\item
there exists $x\in\PP^1\moins U$ which is either a special point or an apparent singularity.
\end{itemize}
There exists a finite number of automorphisms of $\PP^1$ sending $3$ points of $\PP^1\moins U$ to $0,1,\infty$.
\begin{itemize}
\item
We can thus (and will) assume that $0,1,\infty\notin U$ and $\infty$ is either the special point~or, if~such a point does not exist, an apparent singularity.
\end{itemize}

It is then enough to show the finiteness
\begin{enumeratea}
\item\label{enum:a}
of the set of quasi-unipotent rigid irreducible $(\cV,\nabla)$ of rank $r$ with data $(U,\Phi,N)$ and having no special point and an apparent singularity at $\infty$,
\item\label{enum:b}
and of the set of quasi-unipotent rigid irreducible $(\cV,\nabla)$ of rank $r$ having data $(U,\Phi,N)$ such that $\infty$ is special.
\end{enumeratea}

\begin{proof}[Proof of \eqref{enum:a}]
Let us recall the A-D-K algorithm in this case. One shows that there exists a rank-one algebraic bundle with connection $(\cL,\nabla)$ on $U$, completely determined by $(\cV,\nabla)$,
\begin{itemize}
\item
satisfying $\eqref{th:finiteness}_1(U,-\Phi,N)$,
\item
and having monodromy $\chi$ at $\infty$, for some $\chi\in\mu_N\moins\{1\}$, so that $\cV\otimes\cL$ has a regular singularity with monodromy $\chi\id$ at $\infty$,
\end{itemize}
such that $\MC_\chi(\cV\otimes\cL,\nabla)$ has rank $<r$, where $\MC_\chi$ denotes the middle convolution functor by the Kummer sheaf $\cK_\chi$. We note that $\cV\otimes\cL$ has data $(U,\Phi',N)$ with $\Phi'=\Phi-\Phi$.

The formulas given in \cite[Prop.\,1.3.11]{D-S12a} show the following:
\begin{itemize}
\item
$\MC_\chi(\cV\otimes\cL,\nabla)$ has singularities contained in $\PP^1\moins U$, the singularity at $\infty$ being regular with monodromy $\chi^{-1}\id$,
\item
at each $x\in\PP^1\moins(U\cup\{\infty\})$, the set of irregular values of $\MC_\chi(\cV\otimes\cL,\nabla)$ at $x$ is equal to that of $\cV\otimes\cL$ at $x$,
\item
the eigenvalues of the formal monodromies at any $x\in\PP^1\moins U$ belong to $\mu_N$.
\end{itemize}
It follows that $\MC_\chi(\cV\otimes\cL,\nabla)$ has data $(U,\Phi',N)$. Since $\cV\otimes\cL$ is non constant because it is irreducible and has rank $\geq2$, it can be recovered as \hbox{$\MC_{\chi^{-1}}\bigl(\MC_\chi(\cV\otimes\cL,\nabla)\bigr)$}, according to \cite[Prop.\,1.1.9]{D-S12}. We conclude by induction on $r$ that $(\cV,\nabla)$ belongs to the set obtained from a finite set of bundles with connections on~$U$ (satisfying $\eqref{th:finiteness}_{<r}(U,\Phi',N)$) by applying $\MC_{\chi^{-1}}$, for some $\chi\in\mu_N$ and by tensoring by a rank-one bundle with connection belonging to the finite set of those satisfying $\eqref{th:finiteness}_1(U,\Phi,N)$. This shows finiteness in Case \eqref{enum:a}.
\end{proof}

\begin{proof}[Proof of \eqref{enum:b}]
For $(\cV,\nabla)$ having a special point at $\infty$, the A-D-K algorithm starts by exhibiting a summand $\El(\varphi)\otimes R_\varphi$ of $(\cV,\nabla)_{\wh\infty}$ which is ramified of order $\geq2$. It~could occur that, writing $\varphi$ as a minimal ramified polar part $\sum_{j\geq1}^qa_jz^{j/p}$, the leading term $a_qz^{q/p}$ has exponent $q/p$ (the slope of $\varphi$) which is an integer. If such is the case, we~set $\Phi_1=\Phi\cup\{\varphi_1\}$, where $\varphi_1$ is the leading part of $\varphi$ with integral exponents, or~is zero if the slope of $\varphi$ is not an integer. Case II in \cite{Arinkin08} shows that there exists a rank-one algebraic bundle with connection $(\cL,\nabla)$ on $U$, completely determined by $(\cV,\nabla)$, in particular satisfying $\eqref{th:finiteness}_1(U,-\Phi_1,N)$, such that the Fourier transform $\Fou(j_{\dag+}(\cV\otimes\cL,\nabla))$ has rank $<r$ and data $(U',N',\Phi')$. The formal stationary phase formula of \cite[(1.3.5)]{D-S12a} (after \cite{Fang07, Bibi07a, Graham-Squire13}) gives the precise way $(U',\Phi',N')$ is obtained from the data of $(\cV\otimes\cL,\nabla)$ (and thus depends only of these), hence from $(U,\Phi,N)$. One notices that $U'$ and~$\Phi'$ both depend on $(U,\Phi)$ (but not only on $U$ \resp $\Phi$ separately). We conclude by induction on $r$ that $(\cV,\nabla)$ belongs to the set obtained from a finite set of bundles with connections on~$U'$ (satisfying $\eqref{th:finiteness}_{<r}(U',\Phi',N')$) by applying inverse Fourier transform and twist by a finite set of rank-one algebraic bundles with connection, ending thereby the proof for Case \eqref{enum:b}.
\end{proof}

\section{Exponential-geometric origin}
In this section, we prove Property \ref{th:expmot}. Let $\cM$ be the the minimal extension of $(\cV,\nabla)$ on $\PP^1$. It is a quasi-unipotent rigid holonomic $\cD_{\PP^1}$\nobreakdash-mod\-ule. Let us recall the basic result that we will use for the proofs of Theorems \ref{th:expmot} and \ref{th:integrality}. As the proof in \cite{Bibi15} is written in a sketchy way, we give a detailed proof in the appendix.

\begin{proposition}[{\cite[Prop.\,2.69]{Bibi15}}]\label{prop:rigidalgo}
Let $\cM$ be a quasi-unipotent rigid holonomic $\cD_{\PP^1}$\nobreakdash-mod\-ule. There exist
\begin{enumeratea}
\item\label{prop:rigidalgoa}
a smooth projective complex variety $X$ and a strict normal crossing divisor $D\subset X$, together with a subdivisor $D_1\subset D$,
\item\label{prop:rigidalgob}
a projective morphism $f:X\to\PP^1$,
\item\label{prop:rigidalgoc}
a rational function $g$ on $X$ with poles contained in $D$ and whose pole and zero divisors do not intersect,
\item\label{prop:rigidalgod}
a locally free rank-one $\cO_X(*D)$-module $\cN=\cN_\reg$ with a regular singular meromorphic connection $\nabla$,
\end{enumeratea}
\noindent
such that $\cN$ is of torsion (\ie $\cN^{\otimes N}\simeq(\cO_X(*D),\rd)$ for some $N\geq1$) and $\cM$ is the image of the natural morphism
\begin{starequation}
\label{eq:f!*}
\cH^0f_+(\cN,\nabla+\rd g)(!D_1)\to \cH^0f_+(\cN,\nabla+\rd g).
\end{starequation}%
\end{proposition}

In \loccit, $\cM$ is assumed to be formally unitary and its is shown that $\cN$ can be chosen of the same kind. Applying the same proof, one finds that if $\cM$ is quasi-unipotent, then~$\cN$ can be chosen to be of torsion (see the appendix). It is therefore enough to show that the right-hand side of \eqref{eq:f!*} is of exponential-geometric origin. Setting $Y=X\moins D$, we can regard $(\cN,\nabla)$ as a rank-one vector bundle with connection on $Y$ and, in the right-hand side of \eqref{eq:f!*}, we regard $f$ as a morphism $f:Y\to\PP^1$.

\begin{lemme}\label{lem:expmot}
There exists a finite morphism $\rho:X'\to X$ such that $D'=\rho^{-1}(D)$ is a divisor with normal crossings, $\rho:Y':=(X'\moins D')\to Y$ is finite étale, and $\rho^+(\cN,\nabla)\simeq(\cO_{Y'},\rd)$.
\end{lemme}

\begin{proof}
Let $\cN^\nabla$ be the rank-one local system of horizontal sections of $\cN^\an$ on $Y^\an$. Since~$\cN^\nabla$ is of torsion, the monodromies of $\cN^\nabla$ around the various irreducible components of $D$ are roots of the unity, and there exists, after \cite[Th.\,17]{Kawamata81}, a finite morphism $X''\to X$ with~$X''$ smooth projective and the pullback $D''$ of $D$ being a normal crossing divisor such that the pullback of~$\cN^\nabla$ extends as a rank-one local system on~$X''$, which is thus also of torsion. This local system becomes trivial after pullback~by some finite étale covering $X'$ of $X''$, and the composition $X'\to X''\to X$ is the desired~$\rho$.
\end{proof}

We conclude that, on $Y$, $(\cN,\nabla)$ is a direct summand of $\rho_+(\cO_{Y'},\rd)$, since \hbox{$\rho:Y'\to Y$} is finite étale. Furthermore, we have $\cH^0(f\circ\rho)_+=\cH^0 f_+\circ\rho_+$ and thus the holonomic $\cD_{\PP^1}$-module $\cH^0f_+(\cN,\nabla+\rd g)$ is a direct summand of $\cH^0(f\circ\rho)_+(\cO_{Y'},\rd+\rd(g\circ\rho))$. Restricting to the open subset $U$ of $\PP^1$ where these holonomic $\cD_{\PP^1}$-modules are smooth shows that $(\cV,\nabla)$ is of exponential-geometric origin.\qed

\section{Integral structure}
In this section we prove Property \ref{th:integrality}.

\subsection{Stokes-filtered local systems with a \texorpdfstring{$\ZZ[\mu_N]$}{Z}-structure}

Let us set $\kk=\QQ[\zeta]$ for some complex $N$-th root $\zeta$ of $1$ and some $N\geq1$, and let $\bmo=\ZZ[\zeta]$ denote its ring of integers.

Let $\cL$ be a local system of finite-dimensional $\kk$-vector spaces on $U^\an$. It extends in a unique way as a local system of finite-dimensional $\kk$-vector spaces on the real blow-up space $\wt\PP^1$ of $\PP^1$ at each $x\in\PP^1\moins U$. We still denote by $\cL$ this extended $\kk$-local system. For each $x\in\PP^1\moins U$, let $S^1_x$ denote the fiber at $x$ of the real blowing-up map $\varpi:\wt\PP^1\to\PP^1$, and let $\cL(x)$ be the restriction of $\cL$ to $S^1_x$.

In order to define the notion of a Stokes filtration on each $\cL(x)$, we first recall the notion of order between ramified polar parts in a specific direction.

Let $\Phi\subset\CC\lpr t^{1/p}\rpr/\CC\lcr t^{1/p}\rcr$ be a finite set of ramified polar parts and let
\[
\rho: S^1_{x,p}\to S^1_x,\quad z\mto t=z^p
\]
denote the cyclic covering of order $p$. We consider $\Phi$ as a subset of $\CC\lpr z\rpr/\CC\lcr z\rcr$ and we denote by $\theta\bmod2\pi$ a point on $S^1_{x,p}$. For $\varphi,\psi\in\Phi$ we consider the partial order in the direction $\theta$, where we denote by $\nb(\theta)$ a small open sector $(\theta-\epsilon,\theta+\epsilon)\times(0,\delta)$ considered as an open subset of the punctured disc with coordinate $z$:
\[
\psi\leq_\theta\varphi\iff \psi=\varphi\text{ or }\reel(\psi-\varphi)<0\text{ on }\nb(\theta),
\]
and $\psi<_\theta\varphi$ if $\psi\leq_\theta\varphi$ and $\psi\neq\varphi$. The subsets $\{\psi\leq\varphi\}=\{\theta\in S^1_{x,p}\mid\psi\leq_\theta\varphi\}$ and $\psi<\varphi$ are a union of open intervals in $S^1_{x,p}$ and we denote by $\beta_{\psi\leq\varphi}$ the functor composed of the restriction to $\{\psi\leq\varphi\}$ and extension by zero to $S^1_{x,p}$, also denoted $\Gamma_{\{\psi\leq\varphi\}}$; and $\beta_{\psi<\varphi}$ has a similar meaning.

A \emph{graded $\kk$-Stokes-filtered local system index by $\Phi$} is a $\Phi$-graded $\kk$-local system $L_p=\bigoplus_{\varphi\in\Phi}L_{p,\varphi}$ on $S^1_{x,p}$ equipped with the family of nested subsheaves\footnote{We implicitly add the element $-\infty$ to $\Phi$, which satisfies $-\infty<_\theta\varphi$ for any $\varphi\in\Phi$ and any $\theta$, with $L_{p,-\infty}=0$. In such a way, the set $\{\psi\in\Phi\mid\psi<\varphi\}$ is nonempty for any $\varphi$, and $L_{p,<\varphi}$ is possibly zero on some open set.}
\[
L_{p,\leq\varphi}=\bigoplus_{\psi\in\Phi}\beta_{\psi\leq\varphi}L_{p,\psi},\quad
L_{p,<\varphi}=\bigoplus_{\psi\in\Phi}\beta_{p,\psi<\varphi}L_\psi.
\]
Clearly, the following properties are satisfied:
\begin{itemize}
\item
$L_{p,\leq\varphi}/L_{p,<\varphi}=L_\varphi$,
\item
for each $\theta\in S^1_{x,p}$, the family $(L_{p,\leq\varphi,\theta})_{\varphi\in\Phi}$ is an exhaustive increasing filtration\footnote{By exhaustive we mean that $L_{p,\theta}=\bigcup_{\varphi\in\Phi}L_{p,\leq\varphi,\theta}$ for any $\theta$.} with respect to the partial order $\leq_\theta$,
\item
for each $\theta$, we have
\begin{equation}\label{eq:<varphi}
L_{p,<\varphi,\theta}=\sum_{\psi<_\theta\varphi}L_{p,\leq\psi,\theta}\quad\text{(sum taken in $L_{p,\theta}$)}.
\end{equation}
\end{itemize}

It is harmless to enlarge $\Phi$ by adding some ramified polar part $\eta$ and set \hbox{$L_{p,\eta}=0$}. In~such a way, we can (and will implicitly) assume that $\Phi$ is invariant by the automorphisms induced by $z\mto\nu z$ with $\nu^p=1$.

\begin{definition}\label{def:Stokesfiltr}
A $\kk$-Stokes filtration $\cL(x)_\bbullet$ indexed by $\Phi$ of the local system $\cL(x)$ on~$S^1_x$ consists of a family $(\cL(x)_{\leq\varphi})_{\varphi\in\Phi}$ of subsheaves of $\kk$-vector spaces of the local system $\rho^{-1}\cL(x)$ on $S^1_{x,p}$ such that
\begin{enumerate}
\item\label{def:Stokesfiltr2}
\emph{locally on $S^1_{x,p}$}, the pair $(\cL(x),\cL(x)_\bbullet)$ is isomorphic to that of a graded $\kk$\nobreakdash-Stokes-filtered local system of finite-dimensional vector spaces,
\item\label{def:Stokesfiltr3}
for any automorphism $\sigma:S^1_{x,p}\isom S^1_{x,p}$ induced by $z\mto\nu z$ with $\nu^p=1$, and for any $\varphi\in\Phi$, denoting by $a_\sigma(x)$ the canonical identification $\sigma^{-1}\rho^{-1}\cL(x)\simeq\rho^{-1}\cL(x)$, the two subsheaves $a_\sigma(x)(\sigma^{-1}\cL(x)_{\leq\varphi})$ and $\cL(x)_{\leq\sigma^*\varphi}$ of $\rho^{-1}\cL(x)$ are equal.
\end{enumerate}
\end{definition}

\begin{remarque}
From \ref{def:Stokesfiltr}\eqref{def:Stokesfiltr2} and the properties of a graded $\kk$-Stokes-filtered local system, we deduce that, for each $\theta\in S^1_{x,p}$, the germs $\cL_{\leq\varphi,\theta}$ ($\varphi\in\Phi$) are ordered by inclusion according to the partial order $\leq_\theta$ of their indices.

Furthermore, for each $\varphi\in\Phi$, there exists a subsheaf $\cL_{<\varphi}$ well-defined by a formula analogous to \eqref{eq:<varphi} and $\gr_\varphi\cL:=\cL_{\leq\varphi}/\cL_{<\varphi}$ is a locally constant sheaf on $S^1_{x,p}$. As a consequence, $(\cL(x),\cL(x)_\bbullet)$ is locally isomorphic to the graded $\kk$-Stokes-filtered local system $(\gr\cL(x),\gr\cL(x)_\bbullet)$.

Lastly, Property \ref{def:Stokesfiltr}\eqref{def:Stokesfiltr3} also applies to the subsheaves $\cL_{<\varphi}$.
\end{remarque}

Given a ramification $\rho':z'\mto z=z^{\prime q}$, the pullback of a $\kk$-Stokes-filtered local system $(\cL(x),\cL(x)_\bbullet)$ indexed by $\Phi$ is a $\kk$-Stokes-filtered local system indexed by $\rho^{\prime*}\Phi$ and, conversely, any $\kk$-Stokes-filtered local system indexed by $\rho^{\prime*}\Phi$ which is invariant by the automorphisms induced by $z'\mto \nu'z'$ with $\nu^{\prime q}=1$ comes by pullback of a $\kk$-Stokes-filtered local system indexed by $\Phi$.

Given two $\kk$-Stokes-filtered local systems, we can assume that they are indexed by the same $\Phi$. A morphism of $\kk$-Stokes-filtered local systems $(\cL(x),\cL(x)_\bbullet)\to(\cL'(x),\cL'(x)_\bbullet)$ is then a morphism between the corresponding $\kk$-local systems whose pullback by $\rho$ is compatible with the Stokes filtration, and in particular induces a morphism of the corresponding graded $\kk$-local systems.

These notions can be globalized to $U$: a $\kk$-Stokes-filtered local system on $U$ indexed by $\Phi$ consists of a $\kk$-local system on $U^\an$ together with a Stokes filtration indexed by~$\Phi$ on each $\cL(x)$ for $x\in\PP^1\moins U$. A morphism is defined correspondingly.

\begin{theoreme}\label{th:abelian}
The category of $\kk$-Stokes-filtered local system on $U$ is abelian.
\end{theoreme}

\begin{proof}
Since the category of $\kk$-local system on $U^\an$ is abelian, it is enough to consider the category of $\kk$-Stokes-filtered local system on $S^1_x$ ($x\in\PP^1\moins U$). This is \eg \cite[Th.\,3.1]{Bibi10}.\end{proof}

\begin{definition}\label{def:ostructure}
Let $(\cL,(\cL(x)_\bbullet)_{x\in\PP^1\moins U})$ be a $\kk$-Stokes-filtered local system on~$U$ with exponential factors contained in a finite set $\Phi$ of ramified polar parts. An $\bmo=\ZZ[\zeta]$-structure on $(\cL,(\cL(x)_\bbullet)_{x\in\PP^1\moins U})$ consists of
\begin{enumerate}
\item\label{def:ostructure1}
a local system $\cL_{\bmo}$ of $\bmo$-modules of finite type on $U^\an$ such that $\cL=\kk\otimes_{\bmo}\cL_{\bmo}$,
\item\label{def:ostructure2}
for each $x\in\PP^1\moins U$, each $\varphi\in\Phi$ and each $\sigma$  as in Definition \ref{def:Stokesfiltr}\eqref{def:Stokesfiltr3}, a morphism $\lambda_\varphi(x):\cL_{\bmo}(x)_{\leq\varphi}\to\rho^{-1}\cL_{\bmo}(x)$, where each $\cL_{\bmo}(x)_{\leq\varphi}$ is a sheaf of $\bmo$-modules of finite type, and an isomorphism $a_{\varphi,\sigma}(x):\sigma^{-1}\cL_{\bmo}(x)_{\leq\varphi}\isom\cL_{\bmo}(x)_{\leq\sigma^*\varphi}$ such that
\begin{enumerate}
\item\label{def:ostructure2c}
for all $\varphi,\sigma,\sigma'$ we have, with $a_\sigma(x)$ as in Definition \ref{def:Stokesfiltr}\eqref{def:Stokesfiltr3},
\begin{align*}
\lambda_{\sigma^*\varphi}(x)\circ a_{\varphi,\sigma}(x)&=a_{\sigma}(x)\circ\sigma^{-1}\lambda_\varphi(x),
\\
a_{\sigma^*\varphi,\sigma'}(x)\circ \sigma^{\prime-1}a_{\varphi,\sigma}(x)&=a_{(\sigma\sigma')^* \varphi}(x);
\end{align*}
\item\label{def:ostructure2a}
$\id\otimes\lambda_\varphi(x):\kk\otimes_{\bmo}\cL_{\bmo}(x)_{\leq\varphi}\to\kk\otimes_{\bmo}\cL_{\bmo}|_{S^1_{x,p}}=\cL|_{S^1_{x,p}}$ induces an isomorphism onto $\cL(x)_{\leq\varphi}$.
\end{enumerate}
\end{enumerate}

A morphism between $\kk$-Stokes-filtered local systems with an $\bmo$-structure is a morphism between the corresponding $\bmo$-sheaves compatible with the morphisms $\lambda_\varphi(x)$ and $a_{\varphi,\sigma}(x)$ at each $x\in\PP^1\moins U$.
\end{definition}

\begin{remarque}
In Definition \ref{def:ostructure}, we do not impose that the sheaves $\gr_\varphi\cL_\bmo(x)$ are local systems of $\bmo$-modules on $S^1_{x,p}$. This is why we do not use the terminology ``$\bmo$-Stokes-filtered local system''.
\end{remarque}

\begin{corollaire}
The category of $\kk$-Stokes-filtered local systems with an $\bmo$-structure is abelian.
\end{corollaire}

\begin{proof}
The category consisting of objects $(\cL_\bmo,(\cL_\bmo(x)_\bbullet)_{x\in\PP^1\moins U})$ satisfying \ref{def:ostructure}\eqref{def:ostructure1} and~\eqref{def:ostructure2c} is an abelian category. The condition that it yields via \eqref{def:ostructure2a} a $\kk$-Stokes-filtered local system by tensoring the objects with~$\kk$ does not break abelianity, according to Theorem \ref{th:abelian} and $\bmo$-flatness of $\kk$.
\end{proof}

\begin{remarque}[Extension of scalars]\label{rem:extscal}
One can define similarly the notion of a $\CC$-Stokes-filtered local system with a $\kk$-structure and obtain the corresponding abelian category. We notice that the latter category is equivalent to the abelian category obtained from the category of $\kk$-Stokes-filtered local systems by the extension of scalars from $\kk$ to~$\CC$. Indeed, it is a matter of proving that, for a sheaf $\cF_\kk$ of finite-dimensional $\kk$-vector spaces on a locally path connected topological space~$Z$, $\cF_\kk$ is locally constant if and only if $\cF_\CC=\CC\otimes_\kk\cF_\kk$ is so (this is mostly obvious).

As a consequence, the category of $\CC$-Stokes-filtered local system with an $\bmo$-structure is equivalent to the abelian category obtained from the category of $\kk$-Stokes-filtered local systems with an $\bmo$-structure by the extension of scalars from $\kk$ to $\CC$.
\end{remarque}

\begin{remarque}[$\bmo$-structures and Stokes matrices]
The Stokes matrices (or Stokes multipliers) of a $\CC$-Stokes-filtered local system with a $\kk$-structure, equivalently a $\kk$\nobreakdash-Stokes-filtered local system, are conjugate to matrices having entries in $\kk$. On the other hand, in presence of an $\bmo$-structure, we cannot assert in general the existence of Stokes matrices with entries in $\bmo$.

However, in the case of confluent hypergeometric systems considered in Example~\ref{ex:hypergeom}, the computation of the Stokes matrices after a suitable ramification done in \cite{Hien22} provides Stokes matrices with entries in $\bmo$ if the local formal monodromies belong to $\bmo$.

On the other hand, given a locally constant sheaf of \emph{free} $\bmo$-modules on a punctured~$\PP^1$ (without any assumption of irreducibility or rigidity, but one can add them), the computation of the Stokes matrices of the Fourier transform of its associated perverse sheaf on $\Afu$ done in \cite{D-H-M-S17} also provides Stokes matrices with entries in $\bmo$. Such an example, with $\bmo=\ZZ$, can be obtained as follows. Let $f:Y\to\Afu$ be a regular function on a smooth affine complex variety $Y$ of dimension $n$. Assume that $f$ is cohomologically tame (in the sense of \cite{Bibi96bb}), so that in particular $f$ has only isolated critical points in $Y$. Then the Stokes matrices at $t=\infty$ of the free $\CC[t,t^{-1}]$-module with connection
\[
(\cV,\nabla)=\biggl(\frac{\Omega^n(Y)[t,t^{-1}]}{(\rd_Y+t\rd f)\Omega^{n-1}(Y)[t,t^{-1}]},\rd_t+f\rd t\biggr)
\]
can be defined over $\ZZ$. This result goes back to \cite{Pham85b}. Note that $(\cV,\nabla)$, which is of exponential-geometric origin by definition (with $\varphi=tf$ on $Y\times\Gm$), is known to be semi-simple, but is possibly not rigid.
\end{remarque}

\subsection{Proof of Property \ref{th:integrality}}
Let us consider the data as in Proposition \ref{prop:rigidalgo}. Up to blowing up $X$, we can achieve the following properties. There exist a Zariski dense open subset $U$ of $\PP^1$ and a diagram
\[
\xymatrix{
Y=X\moins D=X_U\moins H_U\subset\hspace*{-.9cm}& X_U\ar@{^{ (}->}[r]\ar[d]_{f_U}\ar@{}[rd]|\square&X\ar[d]^f\\
&U\ar@{^{ (}->}[r]&\PP^1
}
\]
such that
\begin{enumeratei}
\item\label{enum:i}
the strict normal crossing divisor $D$ decomposes as $H\cup P$ with $P=f^{-1}(\PP^1\moins U)$,
\item\label{enum:ii}
the pole divisor $P_g$ of $g$ decomposes correspondingly as $P_g=H'\cup P'$ (and the zero divisor of $g$ does not cut $P_g$),
\item\label{enum:iii}
the pair $(X_U^\an,H_U^\an)$ is smooth over $U^\an$, \ie $f^\an$ is smooth on $X_U^\an$ and its restriction to each stratum of the natural stratification of $H_U^\an$ is smooth.
\end{enumeratei}

\begin{lemme}
Assume that the data $(X,f,D,g)$ of Proposition \ref{prop:rigidalgo} satisfy the properties \eqref{enum:i}--\eqref{enum:iii} above, and let $\cN$ be a torsion locally free $\cO_X(*D)$-module with a regular singular meromorphic connection $\nabla$. Let $D_1$ be a sub-divisor of $D$. Then each term of~\eqref{eq:f!*} is a vector bundle with connection on $U$ whose associated $\CC$-Stokes-filtered local system admits a $\ZZ[\zeta]$-structure for which the morphism associated to that of~\eqref{eq:f!*} is a morphism of $\ZZ[\zeta]$-structures.
\end{lemme}

The idea of the proof is that the morphism between the $\CC$-Stokes-filtered local systems associated to \eqref{eq:f!*} can be computed in a purely topological way, by considering suitable real oriented blow-up spaces, from the local system $\cN^\nabla$. The latter being defined over $\bmo=\ZZ[\zeta]$, it follows that the corresponding Stokes-filtered local systems and the morphism between them have an $\bmo$-structure.

\begin{proof}[Proof, part one: the local systems on $U^\an$]
We start with the local systems on $U^\an$. Let us decompose $D_1=H_1\cup P_1$. We~further decompose $H_1$ as $H'_1\cup H''_1$ such that the components of $H'_1$ are the components of $H_1$ along which~$g$ has a pole.

Denoting by $\pDR$ the analytic de~Rham functor (shifted by the ambient dimension), it suffices to show that $\pDR\eqref{eq:f!*}$ is a morphism of (shifted) local systems on $U^\an$ defined over $\bmo$. Due to the commutation of $\pDR$ and direct images, this morphism reads over $U$:
\[
\pR^0f_{U*}\pDR\bigl[(\cN_U,\nabla+\rd g)(!H_{1,U})\bigr]\to \pR^0f_{U*}\pDR(\cN_U,\nabla+\rd g),
\]
and, setting $\cE^g=(\cO_X(*P_g),\rd+\rd g)$, it also reads
\begin{multline}\label{eq:DRf!*}
\pR^0f_{U*}\pDR\bigl[(\cN_U\otimes\cE_U^g)(!H''_{1,U})\bigr] \simeq
\pR^0f_{U*}\pDR\bigl[\cN_U(!H''_{1,U})\otimes\cE_U^g\bigr]\\
\to \pR^0f_*\pDR(\cN_U\otimes\cE_U^g).
\end{multline}

Since the zero and the pole divisor of $g$ do not intersect, the characteristic varieties of $(\cN_U\otimes\cE_U^g)$ and of $(\cN_U\otimes\cE_U^g)(!H_{1,U})\simeq(\cN_U\otimes\cE_U^g)(!H''_{1,U})$ are union of conormal bundles to the strata of the natural stratification of $(X_U,H_U)$. It follows from Assump\-tion \eqref{enum:iii} that the cohomology sheaves of $\bR f_{U*}\pDR\bigl[(\cN_U,\nabla+\rd g)(!H_{1,U})\bigr]$ and $\bR f_{U*}\pDR(\cN_U,\nabla+\rd g)$ are locally constant of finite rank. As a consequence, these cohomology sheaves are, up to a shift, equal to the perverse cohomology sheaves of these complexes.

We will compute these perverse sheaves on $U^\an$ and the morphism between them by means of the real blowing up $\varpi_U:\wt X_U\to X_U$ of the irreducible components of~$H_U$ and of $\wt f_U=f_U\circ\varpi_U:\wt X_U\to U^\an$. Let $\wt X_U^\rmod(g)$ denote the dense open subset of~$\wt X_U$ consisting of points in the neighborhood of which the function $\rme^{-g}$ has moderate growth.
Let $\alpha:Y^\an\hto\wt X_U^\rmod(g)$ and $\beta:\wt X_U^\rmod(g)\hto\wt X_U$ denote the open inclusions, let $\wt f_U^\rmod:X_U^\rmod(g)\to U^\an$ denote the restriction of $\wt f_U$, and let~us set $n=\dim X$. One can define the moderate de~Rham complex $\pDR^\rmod(\cN_U\otimes\cE_U^g)$ on~$\wt X_U$ and a simple computation shows that it is isomorphic to $\beta_!\alpha_*\cN^\nabla_U[n]$, with $\alpha_*\cN^\nabla_U$ being a locally constant sheaf on $\wt X_U^\rmod$. Furthermore, we have a natural isomorphism
\[
\bR\varpi_*\pDR^\rmod(\cN_U\otimes\cE_U^g)\isom\pDR(\cN_U,\nabla+\rd g).
\]
Therefore,
\begin{align*}
\bR f_{U*}\pDR(\cN_U,\nabla+\rd g)&\simeq \bR f_{U*}\bR\varpi_*\pDR^\rmod(\cN_U\otimes\cE_U^g)\\
&=\bR f_{U!}\bR\varpi_!\pDR^\rmod(\cN_U\otimes\cE_U^g)\\
&=\bR \wt f_{U!}\pDR^\rmod(\cN_U\otimes\cE_U^g)\\
&\simeq\bR \wt f_{U!}\bR\beta_!(\alpha_*\cN^\nabla_U)\simeq\bR \wt f_{U!}^\rmod(\alpha_*\cN^\nabla_U).
\end{align*}
It follows that
\begin{itemize}
\item
the locally constant sheaf $\alpha_*\cN^\nabla_U$ has the $\bmo$-structure $\alpha_*\cN_{\bmo,U}^\nabla$,
\item
the sheaf $\pR^0f_*\pDR(\cN_U\otimes\cE_U^g)[-1]$ is isomorphic to the sheaf
\[
R^{n-1}(\wt f_U^\rmod)_!(\alpha_*\cN^\nabla_U)=\CC\otimes_\bmo R^{n-1}(\wt f_U^\rmod)_!(\alpha_*\cN_{\bmo,U}^\nabla).
\]
\end{itemize}
We claim that $R^{n-1}(\wt f_U^\rmod)_!(\alpha_*\cN_{\bmo,U}^\nabla)$ is a locally constant sheaf of $\bmo$-modules of finite type. For that purpose, it is enough to show that $\wt f_U^\rmod:\wt X^\rmod_U\to U^\an$ is a locally trivial topological fibration. Recall that $\wt X_U$ is a $C^\infty$ manifold with corners. It~is thus enough to show that, locally on $U^\an$, any $C^\infty$ vector field on $U^\an$ can be lifted to a $C^\infty$ vector field on $\wt X_U$ tangent to the corners, and whose flow locally preserves $\wt X_U^\rmod$. By using a $C^\infty$ partition of the unity, this is a local question on~$X_U$. By~\eqref{enum:iii}, we can assume that $f_U$ is the second projection $(\CC^{n-1},0)\times(\CC,0)\to(\CC,0)$ and $H_U=(\{x_1\cdots x_\ell=0\},0)\times(\CC,0)$. Furthermore, in this local chart, $g$ can be expressed as the product of a nowhere vanishing holomorphic function by a monomial $x^{-m}$ with $m_i\geq0$ for $i=1,\dots,n-1$. In such a local setting, the description of $\wt X_U$ and $\wt X_U^\rmod$ can be made explicit and the existence of such a lifting of a $C^\infty$ vector field on $(\CC,0)$ is straightforward.

We compute similarly the left-hand side of \eqref{eq:DRf!*}. We denote by $\wt X_{U,1}^\rmod$ the complement of $\varpi_U^{-1}(H''_{1,U})$ in $\wt X_U^\rmod$. We then have the corresponding open inclusions
\[
\alpha_1:Y^\an\hto\wt X_{U,1}^\rmod\quand\beta_1:\wt X_{U,1}^\rmod\hto\wt X_U,
\]
and a computation as above shows that the sheaf $\pR^0f_*\pDR\bigl[\cN_U(!H''_{1,U})\otimes\cE_U^g\bigr][-1]$ is isomorphic to the sheaf
\[
R^{n-1}(\wt f_{U,1}^\rmod)_!(\alpha_{1,*}\cN^\nabla_U)=\CC\otimes_\bmo R^{n-1}(\wt f_{U,1}^\rmod)_!(\alpha_{1,*}\cN_{\bmo,U}^\nabla).
\]
We argue as above to prove that $R^{n-1}(\wt f_{U,1}^\rmod)_!(\alpha_{1,*}\cN_{\bmo,U}^\nabla)$ is a locally constant sheaf of $\bmo$\nobreakdash-modules of finite type.

Letting $\gamma:\wt X_{U,1}^\rmod\hto \wt X_U^\rmod$ denote the open inclusion, we have a natural morphism
\begin{multline*}
R^{n-1}(\wt f_{U,1}^\rmod)_!(\alpha_{1,*}\cN_{\bmo,U}^\nabla)=R^{n-1}(\wt f_U^\rmod)_!(\gamma_!\alpha_{1,*}\cN_{\bmo,U}^\nabla)\\
\to R^{n-1}(\wt f_U^\rmod)_!(\gamma_*\alpha_{1,*}\cN_{\bmo,U}^\nabla)=R^{n-1}(\wt f_U^\rmod)_!(\alpha_*\cN_{\bmo,U}^\nabla),
\end{multline*}
showing that the morphism \eqref{eq:DRf!*} is also defined over $\bmo$.
\end{proof}

\begin{proof}[Proof, part 2: the local Stokes structures]
Since both terms of \eqref{eq:f!*} are holonomic $\cD_{\PP^1}$-modules, their exponential factors at the points of $\PP^1\moins U$ are contained in a finite set of ramified polar parts $\Phi$, of some ramification order $p$. We will not need to compute explicitly this set and we will only use its existence.

We fix $x\in\PP^1\moins U$ and restrict the setting over a small disc $\Delta$ centered at $x$. We~restrict all the data of the lemma as analytic data over $\Delta$. In particular, $U$ is replaced with the punctured disc $\Delta^*$. Otherwise, we keep the same notation with this new analytic meaning.

Since we are only interested in computing Stokes filtrations, we consider the localized modules associated to both terms of \eqref{eq:f!*}, that we regard as meromorphic flat bundles $(\cV_1,\nabla)$ and $(\cV,\nabla)$, \ie free $\cO_\Delta(*0)$-modules of finite rank with a connection, and the natural morphism $(\cV_1,\nabla)\to(\cV,\nabla)$ between them. They have associated $\CC$\nobreakdash-Stokes-filtered local systems $(\cL_1,\cL_{1,\bbullet})$ and $(\cL,\cL_\bbullet)$ with exponential factors contained in $\Phi$. The local systems are the restriction to $\Delta^*$ of those computed in part one of the proof. In particular, we already know that they have an $\bmo$-structure, that we aim at expanding to the whole Stokes structure. For that purpose, we will give a geometric construction of the corresponding Stokes filtrations by means of the maps analogous to $\alpha,\beta,\alpha_1,\beta_1$.

We consider the oriented real blow-ups $\wt X(P)$ and $\wt X=\wt X(D)$ of $X$ along the components of $P$ and $D$ respectively, so that we have a composition
\[
\varpi=\varpi_P\circ\rho:\wt X\to\wt X(P)\to X
\]
and we extend the map $f$ as a continuous map
\[
\wt f:\wt X\xrightarrow{~\rho~}\wt X(P)\to\wt\Delta,
\]
where $\wt\Delta$ is the oriented real blow-up of~$\Delta$ at the origin, with boundary $\partial\wt\Delta\simeq S^1$ and open inclusion $\iota:\Delta^*=\wt\Delta\moins\partial\wt\Delta\hto\wt\Delta$. We~denote by $\wt X_{\leq0}(g)$ the open subset of $\wt X$ consisting of points in the neighborhood of which $\rme^{-g}$ has moderate growth, and we keep the similar notation as in part one of the proof for the maps $\alpha:X\moins D\hto\wt X_{\leq0}(g)$ and $\beta:\wt X_{\leq0}(g)\hto\wt X$.

As in part one, we decompose $D_1$ as $D'_1\cup D''_1$, so that the polar components of $g$ contained in $D_1$ are those of $D'_1$. We consider the open subset
\[
\wt X_{1,\leq0}(g)=\wt X_{\leq0}(g)\moins\varpi^{-1}(D''_1)\subset \wt X_{\leq0}
\]
and we keep the similar notation for the maps $\alpha_1,\beta_1$. Similarly, $\wt f_1$ \resp $\wt f$ denote the restriction of $\wt f$ to $\wt X_{1,\leq0}(g)$ \resp $\wt X_{\leq0}(g)$. From \cite[Cor.\,4.7.5 \& Lem.\,5.1.6]{Mochizuki10} we obtain:

\begin{lemme}\label{lem:mochi}
There exists a commutative diagram
\[
\xymatrix{
R^{n-1}(\wt f_1)_!(\alpha_{1,*}\cN^\nabla)|_{\partial\wt\Delta}\ar[d]\ar@{}[r]|-{\simeq}&\cL_{1,\leq0}\ar[d]\\
R^{n-1}\wt f_!(\alpha_*\cN^\nabla)|_{\partial\wt\Delta}\ar@{}[r]|-{\simeq}&\cL_{\leq0}
}
\]
where the vertical morphisms are the natural ones and, by means of the horizontal isomorphisms, the inclusion $\cL_{\leq0}\subset\cL$ is the adjunction morphism $R^{n-1}\wt f_!(\alpha_*\cN^\nabla)|_{\partial\wt\Delta}\to \bigl(\iota_*\iota^{-1}R^{n-1}\wt f_!(\alpha_*\cN^\nabla)\bigr)|_{\partial\wt\Delta}$, and similarly for $\cL_1$.\qed

\end{lemme}

Since $\cN^\nabla$ is equipped with the $\bmo$-structure $\cN^\nabla_\bmo$, all terms and morphisms in the lemma acquire a natural $\bmo$-structure compatible with that already obtained for $\cL_1,\cL$ via the adjunction morphism.

In order to obtain,  for any $\varphi\in\Phi\subset\CC\lpr t^{1/p}\rpr/\CC\lcr t^{1/p}\rcr$, the $\bmo$-structure on $\cL_{1,\leq\varphi}$ and $\cL_{\leq\varphi}$, we consider the diagram
\[
\xymatrix{
X_p\ar[d]\ar@/_1.5pc/[dd]_{f_p}\\
X'_p\ar[r]\ar[d]\ar@{}[dr]|-\square&X\ar[d]^f\\
\Delta_p\ar[r]_-{\rho_p}&\Delta
}
\]
where $\rho_p$ is the cyclic ramification of order $p$ and $X_p$ is a resolution of singularities of the pair $(X'_p,P)$. We replace the rational function $g$ on $X$ with $g+f_p^*\varphi$ on $X_p$, the divisor $D$ with its pullback by $X_p\to X$, and $\cN$ with its pullback $\cN_p$ on $X_p$. Then $\cL_{1,\leq\varphi},\cL_{\leq\varphi}$ are obtained by the same procedure as that of Lemma \ref{lem:mochi} with these new data, so that the $\bmo$-structure is obtained in the same way. In this case, the morphisms $\lambda_{\varphi}(x)$ of Definition \ref{def:ostructure}\eqref{def:ostructure2} are the adjunction morphisms, and the compatibility of the morphisms $a_{\varphi,\sigma}(x)$ with products $\sigma'\sigma$ follows from the identification
\[
\sigma^{\prime-1}\wt X_{p,\leq0}(g+\sigma^*\varphi)=\wt X_{p,\leq0}(\sigma^{\prime*}(g+\sigma^*\varphi))=\wt X_{p,\leq0}(g+(\sigma\sigma')^*\varphi).
\]
This concludes the proof of Property \ref{th:integrality}.
\end{proof}

\subsection{Final remarks on the integral structure}
There are other possible approaches to $\bmo$\nobreakdash-struc\-tures. They all mainly rely on the general Riemann-Hilbert correspondence as developed by D'Agnolo and Kashiwara \cite{D-K13} and their subsequent work. As indicated in \cite[\S2]{D-K13}, the theory of \loccit\ can be applied to objects defined on $\bmo$. The point is then to check that the $\RR$-constructible enhanced ind-sheaves associated to $\cN\otimes\cE^\varphi$ and $\cN(!D_1)\otimes\cE^\varphi$ and the morphism between them are defined over $\bmo$. The pushforward of these objects by the projective morphism $f$ provides an $\bmo$-structure on the enhanced ind-sheaves associated to both terms of \eqref{eq:f!*} and the morphism between them. This can be made a little more precise by considering the categories of $\CC$-constructible enhanced ind-sheaves of \cite{Ito20}, or the characterization of $\RR$-constructible enhanced ind-sheaves which come from holonomic $\cD$-modules given in \cite{Mochizuki16}.

One could also work within the setting of irregular constructible complexes of \cite{Kuwagaki18}. Such objects can be defined over the ring $\bmo$.

In all these theories, the main point is the compatibility of the irregular Riemann-Hilbert correspondence with projective pushforward. This is probably one of the most delicate points in \cite{D-K13}, that replaces \cite[Cor.\,4.7.5 \& Lem.\,5.1.6]{Mochizuki10} used in Lemma \ref{lem:mochi}.

\appendix\refstepcounter{section}
\section*{Appendix: detailed proof of Proposition \ref{prop:rigidalgo}}
\let\map f
\let\mero g
\let\modif e
\let\regg\reg

\begin{notation*}
Let $H$ be a hypersurface in a smooth variety $X$. We denote by $\Gamma_{[*H]}$ the localization functor acting on the category of holonomic $\cD_X$-modules: for such a $\cD_X$-module $M$, we have $\Gamma_{[*H]}M=\cO_X(*H)\otimes_{\cO_X}M$ as an $\cO_X$-module. Denoting by $\bD$ the duality functor on holonomic $\cD_X$-modules, we define the dual localization functor $\Gamma_{[!H]}$ as $\bD\Gamma_{[*H]}\bD$. Both functors vanish when applied to holonomic $\cD_X$-modules supported on $H$ (this is clear for the localization functor, and the property for the dual one follows from the fact that duality preserves the support).

Let $\mero$ be a meromorphic (or rational) function on $X$. If $M$ is considered as an $\cO_X$-module with an integrable connection $\nabla$, it will be convenient to denote the same $\cO_X$-module with the twisted connection $\nabla+\rd\mero$ as $\cE^\mero\otimes M$. Note that this twist contains the localization functor along the pole divisor $P$ of $\mero$, so that $\cE^\mero\otimes M=\Gamma_{[*P]}(\cE^\mero\otimes M)$.

We will make use of the relations between various functors described in \cite[\S1.77]{Bibi15}.
\end{notation*}

Let $\cM$ be rigid irreducible on $\PP^1$ and let us assume that there exist $X,D,D_1$, $\map:X\to\PP^1$, $\mero$ and $(\cN,\nabla)$ as in \ref{prop:rigidalgo}\eqref{prop:rigidalgoa}\nobreakdash--\eqref{prop:rigidalgod} so that $\cM$ is the image of \eqref{eq:f!*}. Denoting by $P$ the pole divisor of $\mero$, we can assume that $D_1$ has no component contained in $P$: indeed, denoting by $D_0$ the union of those components not contained in $P$, we have
\[
\cE^\mero\otimes\Gamma_{[!D_1]}\cN=\cE^\mero\otimes\Gamma_{[!D_0]}\cN\quand\cE^\mero\otimes\Gamma_{[*D_1]}\cN=\cE^\mero\otimes\Gamma_{[*D_0]}\cN.
\]
We will prove the following properties.
\begin{enumeratei}
\item\label{enum:rigidalgo1}
If $\cL$ is a rank-one meromorphic connection on $\PP^1$ with poles along $\Sigma\subset\PP^1$, then there exist data \ref{prop:rigidalgo}\eqref{prop:rigidalgoa}--\eqref{prop:rigidalgod} such that the image $\cM'$ of
\[
\Gamma_{[!\Sigma]}(\cM\otimes\cL)\to\Gamma_{[*\Sigma]}(\cM\otimes\cL)
\]
is the image of \eqref{eq:f!*} with these data.

Moreover, if $\cN^\nabla$ is unitary (\resp of torsion) and $\cL$ is locally formally unitary (\resp formally quasi-unipotent), then so is $\cM'$ and $(\cN',\nabla)$ can be chosen unitary (\resp of torsion).

\item\label{enum:rigidalgo2}
Let $\Afu_t$ be the chart with coordinate $t$ corresponding to the choice $0,1,\infty\in\PP^1$. By \eqref{enum:rigidalgo1}, we can assume that $\cM=\Gamma_{[*\infty]}\cM$. Let $\cM'$ be the Laplace transform of $\cM$ with respect to this choice. Then there exist data \ref{prop:rigidalgo}\eqref{prop:rigidalgoa}--\eqref{prop:rigidalgod} such that $\cM'$ is the image of \eqref{eq:f!*} with these data.

Moreover, if $\cN^\nabla$ is unitary (\resp of torsion) and $\cL$ is locally formally unitary (\resp formally quasi-unipotent), then so is $\cM'$ and $(\cN',\nabla)$ can be chosen unitary (\resp of torsion).
\end{enumeratei}

These two properties allow us to conclude the proof of the proposition, since any rigid holonomic $\cD_{\PP^1}$-module can be obtained by applying a sequence of \eqref{enum:rigidalgo1} and \eqref{enum:rigidalgo2} to $(\cO_{\PP^1},\rd)$, according to the Arinkin-Deligne algorithm, and moreover, if $\cM$ is locally formally unitary (\resp quasi-unipotent), then the rank-one connections $\cL_\regg$ chosen at each step are locally unitary (\resp quasi-unipotent).

Let us show \eqref{enum:rigidalgo1}. There exists a meromorphic function $\psi$ on $\PP^1$ and a rank-one meromorphic connection $\cL_\regg$ with regular singularities, such that $\cL=\cE^\psi\otimes\cL_\regg$. We can write $\cL=(\cO_{\PP^1}(*\Sigma),\rd+\rd\psi+\omega)$, where $\Sigma$ is the pole divisor of $\cL$ and~$\omega$ is a one-form with at most simple poles at $\Sigma$. Moreover, $\cL$ is locally formally unitary (\resp quasi-unipotent) if and only if $\cL_\regg$ is unitary (\resp quasi-unipotent), \ie the residues of $\omega$ at $\Sigma$ are real (\resp rational).

Since $\cL$ is $\cO_{\PP^1}$-flat, $\cM\otimes\cL$ is the image of
\[
\map^0_\dag(\cE^\mero\otimes\Gamma_{[!D_1]}\cN)\otimes\cL\to \map^0_\dag(\cE^\mero\otimes\Gamma_{[*D_1]}\cN)\otimes\cL,
\]
and, since the functors $\Gamma_{[\star\Sigma]}$ ($\star=!,*$) are exact on the category of holonomic $\cD_{\PP^1}$\nobreakdash-mod\-ules, $\cM'$ is the image of
\[
\Gamma_{[!\Sigma]}\map^0_\dag(\cE^\mero\otimes\Gamma_{[!D_1]}\cN)\otimes\cL\to \Gamma_{[*\Sigma]}\map^0_\dag(\cE^\mero\otimes\Gamma_{[*D_1]}\cN)\otimes\cL.
\]

We set $H=\map^{-1}(\Sigma)$ and we decompose $D$ as $D_2\cup D_3\cup D_4$, where $D_3$ are those components of $D$ which are components of $H$, $D_2$ are the components of $D_1$ which are not components of $H$, and $D_4$ are the remaining components. We set $\mero_1=\mero+\map\circ\psi$.

Firstly, one checks that, for $\star=*,!$,
\[
\map_\dag^0(\cE^\mero\otimes\Gamma_{[\star D_1]}\cN)\otimes\cL\simeq \map_\dag^0\bigl( (\cE^\mero\otimes\Gamma_{[\star D_1]}\cN)\otimes\map^+\cL\bigr)\simeq\map_\dag^0\bigl(\cE^{\mero_1}\otimes\Gamma_{[\star D_1]}\cN\otimes\map^+\cL_\regg\bigr).
\]
Due to the commutation $\Gamma_{[\star\Sigma]}\map^0_\dag\simeq\map^0_\dag\Gamma_{[\star H]}$ for $\star=*,!$, and since $\map^+\cL_\regg=\Gamma_{[*H]}\map^+\cL_\regg$ (so that we can replace $D_1$ with $D_2$ in the right-hand side), we deduce that $\cM'$ is the image~of
\begin{equation}\label{eq:*}
\map_\dag^0\bigl(\Gamma_{[!H]}(\cE^{\mero_1}\otimes\Gamma_{[!D_2]}\cN\otimes\map^+\cL_\regg)\bigr)\to\map_\dag^0\bigl(\cE^{\mero_1}\otimes\cN\otimes\map^+\cL_\regg\bigr).
\end{equation}

Let $\modif:X'\to X$ be a projective modification such that $e^{-1}(D\cup H)$ is a divisor with normal crossings and such that the pole and zero divisors of $\mero'_1:=\modif^*\mero_1$ do not intersect. For the first condition the blowing ups can be chosen to take place above the union of $D\cap H$ and of the singular set of $H$, while for the second condition, since the pole and zero divisors of $\mero_1$ intersect at most in $D\cap H$, the blowing ups can be chosen to take place above $D\cap H$. As a consequence, we can assume that, setting $H'=\modif^{-1}(H)$, the morphism $\modif:X'\moins H'\to X\moins H$ is an isomorphism.

Since $\cE^{\mero_1}\otimes\Gamma_{[!D_2]}\cN\otimes\map^+\cL_\regg=\Gamma_{[*H]}(\cE^{\mero_1}\otimes\Gamma_{[!D_2]}\cN\otimes\map^+\cL_\regg)$, we have, after \cite[1.77(vi)]{Bibi15},
\[
\modif^k_\dag(\cE^{\mero'_1}\otimes\modif^*(\Gamma_{[!D_2]}\cN\otimes\map^+\cL_\regg))\simeq
\begin{cases}
\cE^{\mero_1}\otimes\Gamma_{[!D_2]}\cN\otimes\map^+\cL_\regg&\text{if }k=0,\\
0&\text{otherwise}.
\end{cases}
\]
With \cite[1.77(iii) \& Cor.\,1.80]{Bibi15} we deduce
\[
\modif^k_\dag(\cE^{\mero'_1}\otimes\Gamma_{[!H']}\modif^+(\Gamma_{[!D_2]}\cN\otimes\map^+\cL_\regg))\simeq\begin{cases}
\Gamma_{[!H]}(\cE^{\mero_1}\otimes\Gamma_{[!D_2]}\cN\otimes\map^+\cL_\regg)&\text{if }k=0,\\
0&\text{otherwise}.
\end{cases}
\]
Setting $\map'=\map\circ\modif$, we conclude that the left-hand side of \eqref{eq:*} is isomorphic to
\[
\map_\dag^{\prime0}(\cE^{\mero'_1}\otimes\Gamma_{[!H']}\modif^+(\Gamma_{[!D_2]}\cN\otimes\map^+\cL_\regg)).
\]
On the other hand, the right-hand side of \eqref{eq:*} is easily computed as
\[
\map_\dag^{\prime0}(\cE^{\mero'_1}\otimes\modif^+(\cN\otimes\map^+\cL_\regg)),
\]
and the morphism in \eqref{eq:*} is the natural one with the expressions above.

Let $D'_2$ be the strict transform of $D_2$ by $\modif$. We claim that
\begin{equation}\label{eq:**}
\Gamma_{[!H']}\modif^+(\Gamma_{[!D_2]}\cN\otimes\map^+\cL_\regg)\simeq\Gamma_{[!H']}\Gamma_{[!D'_2]}\modif^+(\cN\otimes\map^+\cL_\regg),
\end{equation}
so that, setting $D_1'=H'\cup D'_2$, we find that $\cM'$ is the image of the natural morphism
\[
\map_\dag^{\prime0}(\cE^{\mero'_1}\otimes\Gamma_{[!D_1']}\modif^+(\cN\otimes\map^+\cL_\regg))\to
\map_\dag^{\prime0}(\cE^{\mero'_1}\otimes\modif^+(\cN\otimes\map^+\cL_\regg)),
\]
and this concludes the proof of \eqref{enum:rigidalgo1}.

Let us prove that \eqref{eq:**} holds true. On the one hand, the kernel and cokernel of the morphism
\[
\modif^+(\Gamma_{[!D_2]}\cN\otimes\map^+\cL_\regg)\to\modif^+(\cN\otimes\map^+\cL_\regg)
\]
are supported on $D'_2$, hence they vanish after applying $\Gamma_{[!D'_2]}$, so that
\[
\Gamma_{[!D'_2]}\modif^+(\Gamma_{[!D_2]}\cN\otimes\map^+\cL_\regg)\to\Gamma_{[!D'_2]}\modif^+(\cN\otimes\map^+\cL_\regg)
\]
is an isomorphism. On the other hand, the natural morphism
\[
\Gamma_{[!D'_2]}\modif^+(\Gamma_{[!D_2]}\cN\otimes\map^+\cL_\regg)\to\modif^+(\Gamma_{[!D_2]}\cN\otimes\map^+\cL_\regg)
\]
is an isomorphism away from $H'$: since $\map^+\cL_\regg$ is locally isomorphic to $(\cO_{X'},\rd)$ on $X'\moins H'$, $\Gamma_{[!D_2]}\cN\otimes\map^+\cL_\regg$ is locally isomorphic to $\Gamma_{[!D_2]}(\cN\otimes\map^+\cL_\regg)$ on this open set, so that the assertion reduces to the equality $\Gamma_{[!D_2]}(\cbbullet)=\Gamma_{[!D_2]}\Gamma_{[!D_2]}(\cbbullet)$. As a consequence, the kernel and cokernel of the above morphism have support in $H'$. Therefore, they vanish after the application of the exact functor $\Gamma_{[!H']}$.\qed

\medskip
Let us now show \eqref{enum:rigidalgo2}. Recall that the Laplace transformation is an exact functor $\Mod_\hol(\PP^1,*\infty)\mto\Mod_\hol(\PP^1,*\infty)$, so $\cM'$ is the image of the Laplace transform of $\map^0_\dag(\cE^\mero\otimes\Gamma_{[!D_1]}\cN)$ to that of $\map^0_\dag(\cE^\mero\otimes\Gamma_{[*D_1]}\cN)$. We denote by $t$ the variable on $\PP^1_t$ and we introduce a new $\PP^1$ with variable $\tau$, and we denote by $p,q:\PP^1_t\times\PP^1_\tau\to\PP^1_t,\PP^1_\tau$ the first and second projections. For any holonomic $\cD_{\PP^1_t}$-module $M$, its Laplace transform is given by the formula $q^0_\dag(\cE^{-t\tau}\otimes\Gamma_{[*(\PP^1_t\times\infty)]}p^+\Gamma_{[*\infty]}M)$ (and $q_\dag^k(\cdots)=0$ for $k\neq0$). We apply this formula to $\cM=\map^0_\dag(\cE^\mero\otimes\Gamma_{[\star D_1]}\cN)$ ($\star=*,!$).

Since we have assumed $\cM=\Gamma_{[*\infty]}\cM$, we can assume that $D$ contains $\map^{-1}(\infty)$ but~$D_1$ does not contain any component of it, and we can omit $\Gamma_{[*\infty]}$ in the previous formula. We denote by $\wh p,\wh q$ the projections $\wh X:=X\times\PP^1_\tau\to X,\PP^1_\tau$, so that we have $q\circ(\map\times\id_{\PP^1_\tau})=\wh q$. For $\star=*,!$, we~have
\begin{align*}
q^0_\dag\Bigl(&\cE^{-t\tau}\otimes \Gamma_{[*(\PP^1_t\times\infty)]}p^+\map^0_\dag(\cE^\mero\otimes\Gamma_{[\star D_1]}\cN)\Bigr)\\
&\simeq q^0_\dag\Bigl(\cE^{-t\tau}\otimes\Gamma_{[*(\PP^1_t\times\infty)]}(\map\times\id)^0_\dag \wh p^+(\cE^\mero\otimes\Gamma_{[\star D_1]}\cN)\Bigr)\quad(\text{after \cite[1.77(v)]{Bibi15}})\\
&\simeq q^0_\dag\Bigl(\cE^{-t\tau}\otimes (\map\times\id)^0_\dag \Gamma_{[*(X\times\infty)]}\wh p^+(\cE^\mero\otimes\Gamma_{[\star D_1]}\cN)\Bigr)\quad(\text{after \cite[1.77(iii)]{Bibi15}})\\
&\simeq q^0_\dag\Bigl((\map\times\id)^0_\dag \Gamma_{[*(X\times\infty)]}(\cE^{p^*\mero-\tau\map}\otimes \wh p^+\Gamma_{[\star D_1]}\cN)\Bigr)\quad(\text{after \cite[1.77(iii)]{Bibi15}})\\
&\simeq \wh q^0_{\dag}\Bigl(\cE^{p^*\mero-\tau\map}\otimes\Gamma_{[*(X\times\infty)]}\wh p^+\Gamma_{[\star D_1]}\cN\Bigr)\quad(\text{since }q^k_\dag(\cdots)=0\text{ for }k\neq0).
\end{align*}
Set $\wh D=(D\times\PP^1_\tau)\cup(X\times\infty)$, $\wh D_1=D_1\times\PP^1_\tau$ and $\wh\cN=\Gamma_{[*(X\times\infty)]}\wh p^+\cN$. We have $\wh \cN=\Gamma_{[*\wh D]}\wh\cN$. Then the latter expression can also be written as
\[
\wh q^0_\dag\bigl(\cE^{p^*\mero-\tau\map}\otimes\Gamma_{[\star \wh D_1]}\wh\cN\bigr),
\]
and therefore $\cM'$ is the image of the natural morphism
\[
\wh q^0_\dag\bigl(\cE^{p^*\mero-\tau\map}\otimes\Gamma_{[! \wh D_1]}\wh\cN\bigr)\to\wh q^0_\dag\bigl(\cE^{p^*\mero-\tau\map}\otimes\Gamma_{[* \wh D_1]}\wh\cN\bigr).
\]

Let $\modif':\wh X'\to \wh X$ be a projective modification which induces an isomorphism above the complement of $\wh D$, whose pullback is denoted by~$\wh D'$, such that the pole and zero divisors of $\mero':=(p\circ \modif')^*\mero-\modif^{\prime*}(\tau\map)$ do not intersect, and set~$\wh D_1'=\modif^{\prime-1}(\wh D_1)$ and $\wh \cN'=\Gamma_{[*\wh D']}e^{\prime+}\wh\cN$. Then, arguing as in Case \eqref{enum:rigidalgo1}, we find
\begin{align*}
\wh q^0_\dag\bigl(\cE^{p^*\mero-\tau\map}\otimes\Gamma_{[\star \wh D_1]}\wh\cN\bigr)
&\simeq \wh q^0_\dag\modif^{\prime0}_\dag\bigl(\cE^{\mero'}\otimes\Gamma_{[\star \wh D_1']}\cN'\bigr)\\
&\simeq (\wh q\circ \modif')^0_\dag\bigl(\cE^{\mero'}\otimes\Gamma_{[\star \wh D_1']}\cN'\bigr)\quad(\text{since }\modif^{\prime k}_\dag(\cdots)=0\text{ for }k\neq0).
\end{align*}
Set now $\map'=\wh q\circ \modif'$. Then, $\cM'$ is the image of \eqref{eq:f!*} with respect to the data $\wh X'$, $\wh D'$, $\wh D_1'$, $\map'$ and $\mero',\wh\cN'$. Furthermore, $\wh\cN'$ is unitary (\resp of torsion) if $\cN$ is so.\qed

\subsubsection*{Acknowledgements}
The author thanks the referees for their careful reading of the manuscript, for noticing various inaccuracies in the first version, and for suggesting various improvements.

\backmatter
\providecommand{\didotfam}{}
\providecommand{\bysame}{\leavevmode ---\ }
\providecommand{\og}{``}
\providecommand{\fg}{''}
\providecommand{\cdrandname}{\&}
\providecommand{\cdredsname}{\'eds.}
\providecommand{\cdredname}{\'ed.}
\providecommand{\cdrmastersthesisname}{M\'emoire}
\providecommand{\cdrphdthesisname}{Th\`ese}
\providecommand{\eprint}[1]{\href{http://arxiv.org/abs/#1}{\texttt{arXiv\string:\allowbreak#1}}}
\providecommand{\eprintother}[3]{\href{#1/#2}{\texttt{#3\string:\allowbreak#2}}}

\end{document}